%% file: nphard.tex
\documentclass[12pt]{article}

\usepackage{amsfonts}
\usepackage{amssymb}
\usepackage{epsfig, subfigure, verbatim}

\renewcommand{\emph}[1]{{\it #1}}
\newcounter{segcount}

\newenvironment{segment}
{\refstepcounter{segcount}\vspace{5mm}
\noindent{\bf \thesegcount. }}
{}

\newenvironment{statementnumbered}[3][]
{\refstepcounter{segcount}\vspace{5mm}
\noindent{\bf\thesegcount. #2}#1{\bf.}\ {\sl #3}}
{\nolinebreak[4] \nopagebreak[4] $\hfill \square$}

\newenvironment{statement}[3][]
{\vspace{5mm}\noindent{\bf #2}#1{\bf.} {\sl #3}}
{\nolinebreak[4] \nopagebreak[4] $\hfill \square$}

\newenvironment{statementnoboxnumbered}[3][]
{\refstepcounter{segcount}\vspace{5mm}
\noindent{\bf\thesegcount. #2}#1{\bf.} {\sl #3}}
{}

\newenvironment{statementnobox}[3][]
{\vspace{5mm}\noindent{\bf #2}#1{\bf.} {\sl #3}}
{}

\newenvironment{definitionnumbered}[2][]
{\refstepcounter{segcount}\vspace{5mm}
\noindent{\bf\thesegcount. #2}#1{\bf.}}
{}

\newenvironment{definition}[2][]
{\vspace{5mm}\noindent{\bf #2}#1{\bf.}}
{}

\newenvironment{resultnumbered}[3][]
{\refstepcounter{segcount}\vspace{5mm}
\noindent{\bf\thesegcount. #2}#1{\bf.} {\sl #3}
\vskip5mm\noindent {\bf Proof: }}
{\nopagebreak[4] $\hfill \square$}

\newenvironment{result}[3][]
{\vspace{5mm}
\noindent{\bf#2}#1{\bf.} {\sl #3}
{\\ \bf Proof: }}
{\nolinebreak[4] \nopagebreak[4] $\hfill \square$}

\newenvironment{risultnumbered}[3][]
{\refstepcounter{segcount}\vspace{5mm}
\noindent{\bf\thesegcount. #2}#1{\bf.} {\sl #3}
{\nopagebreak[4] \noindent \bf Proof: }}
{$\hfill \square$}

\newenvironment{risult}[3][]
{\vspace{5mm}
\noindent{\bf#2}#1{\bf.} {\sl #3}
{\noindent \bf Proof: }}
{\nolinebreak[4] \nopagebreak[4] $\hfill \square$}

\newenvironment{risultnoboxnumbered}[3][]
{\refstepcounter{segcount}\vspace{5mm}
\noindent{\bf\thesegcount. #2}#1{\bf.} {\sl #3}
{\nopagebreak[4] \noindent \bf Proof: }}
{}

\newenvironment{risultnobox}[3][]
{\vspace{5mm}
\noindent{\bf#2}#1{\bf.} {\sl #3}
{\nopagebreak[4] \noindent \bf Proof: }}
{}

\newenvironment{resultnoboxnumbered}[3][]
{\refstepcounter{segcount}\vspace{5mm}
\noindent{\bf\thesegcount. #2}#1{\bf.} {\sl #3}
{\\ \bf Proof: }}
{}

\newenvironment{resultnobox}[3][]
{\vspace{5mm}\noindent{\bf #2}#1{\bf.} {\sl #3}
{\\ \bf Proof: }}
{}

\newcommand{\seg}{\begin{segment}}
\newcommand{\segend}{\end{segment}}

\newcommand{\stmtnum}{\begin{statementnumbered}}
\newcommand{\stmtnumend}{\end{statementnumbered}}

\newcommand{\stmt}{\begin{statement}}
\newcommand{\stmtend}{\end{statement}}

\newcommand{\stmtnoboxnum}{\begin{statementnoboxnumbered}}
\newcommand{\stmtnoboxnumend}{\end{statementnoboxnumbered}}

\newcommand{\stmtnobox}{\begin{statementnobox}}
\newcommand{\stmtnoboxend}{\end{statementnobox}}

\newcommand{\defnnum}{\begin{definitionnumbered}}
\newcommand{\defnnumend}{\end{definitionnumbered}}

\newcommand{\defn}{\begin{definition}}
\newcommand{\defnend}{\end{definition}}

\newcommand{\resnum}{\begin{resultnumbered}}
\newcommand{\resnumend}{\end{resultnumbered}}

\newcommand{\res}{\begin{result}}
\newcommand{\resend}{\end{result}}

\newcommand{\risnum}{\begin{risultnumbered}}
\newcommand{\risnumend}{\end{risultnumbered}}

\newcommand{\ris}{\begin{risult}}
\newcommand{\risend}{\end{risult}}

\newcommand{\risnoboxnum}{\begin{risultnoboxnumbered}}
\newcommand{\risnoboxnumend}{\end{risultnoboxnumbered}}

\newcommand{\risnobox}{\begin{risultnobox}}
\newcommand{\risnoboxend}{\end{risultnobox}}

\newcommand{\resnoboxnum}{\begin{resultnoboxnumbered}}
\newcommand{\resnoboxnumend}{\end{resultnoboxnumbered}}

\newcommand{\resnobox}{\begin{resultnobox}}
\newcommand{\resnoboxend}{\end{resultnobox}}

\newcommand{\cA}{{\cal A}}
\newcommand{\cB}{{\cal B}}
\newcommand{\cD}{{\cal D}}
\newcommand{\cF}{{\cal F}}
\newcommand{\cK}{{\cal K}}
\newcommand{\cO}{{\cal O}}
\newcommand{\cP}{{\cal P}}
\newcommand{\cQ}{{\cal Q}}
\newcommand{\cR}{{\cal R}}

\title{Vertex-partitioning into fixed additive induced-hereditary properties is NP-hard}

\author{
Alastair Farrugia
\\ {\it afarrugia@math.uwaterloo.ca}
\\ Dept. of Combinatorics \& Optimization
\\ University of Waterloo, Ontario, Canada, N2L 3G1
}

\date{December 19, 2002}

\begin{document}
\maketitle
\begin{abstract}
Can the vertices of a graph $G$ be partitioned into $A \cup B$, so that $G[A]$ is a line-graph and $G[B]$ is a forest? Can $G$ be partitioned into a planar graph and a perfect graph? The NP-completeness of these problems are just special cases of our result:
if $\cP$ and $\cQ$ are additive induced-hereditary graph properties, then $(\cP, \cQ)$-colouring is NP-hard, with the sole exception of graph $2$-colouring 
(the case where both $\cal P$ and $\cal Q$ are the set $\cO$ of finite edgeless graphs). Moreover, $(\cP, \cQ)$-colouring is NP-complete iff $\cP$- and $\cQ$-recognition are both in NP. This proves a conjecture of Kratochv\'{\i}l and Schiermeyer.
\end{abstract}

%

Kratochv\'{\i}l and Schiermeyer conjectured in~\cite{comp3} that for any additive hereditary graph properties $\cP$ and $\cQ$, recognising graphs in $\cP\circ \cQ$ is 
NP-hard, with the obvious exception of bipartite graphs (the case where both $\cal P$ and $\cal Q$ are the set $\cal O$ of finite edgeless graphs). They settled the case where $\cQ = \cO$, and it was natural to extend the conjecture to 
\emph{induced}-hereditary properties. 
Berger's result~\cite{berger} that reducible additive induced-hereditary properties have infinitely many minimal forbidden subgraphs provided support for the extended conjecture.

We prove the extension of the Kratochv\'{\i}l-Schiermeyer conjecture in this paper. 
Problems such as the following (for an arbitrary graph $G$) are therefore NP-complete.
Can $V(G)$ be partitioned into $A \cup B$, so that $G[A]$ is a line-graph and $G[B]$ is a forest? Can $G$ be partitioned into a planar graph and a perfect graph? For fixed $k, \ell, m$, can $G$ be partitioned into a $k$-degenerate subgraph, a subgraph of maximum degree $\ell$, and an $m$-edge-colourable subgraph?
\newline

Garey et al.\ \cite{g-j-stock, monien} essentially showed
$(\cO, \{\textrm{forests}\})$-colouring to be NP-complete,
while Brandst\"{a}dt et al.~\cite[Thm. 3]{brand} proved the case 
$(\cO, \{P_4, C_4\}-\textrm{free graphs})$.

Let $\cal P$ be a property and let ${\cal P}^k$ be the product of $\cal P$ with itself, $k$ times. Brown and Corneil~\cite{brown-thesis, brown-perfect} showed that ${\cal P}^k$-recognition is NP-hard when $\cal P$ is the set of perfect graphs and $k\ge 2$, while Hakimi and Schmeichel~\cite{arbor} did the case $\{$forests$\}^2$. 
There was particular interest in $G$-free $k$-colouring (where $\cP$ has just one forbidden induced-subgraph $G$). When $G = K_2$ we get graph colouring, one of the best known NP-complete problems, while subchromatic number~\cite{subchrom, subchrom-2} (partitioning into subgraphs whose components are all cliques) is the case $G = P_3$.
Brown~\cite{brown} proved the case where $G$ is  $2$-connected, and Achlioptas~\cite{comp1} showed NP-completeness for all $G$.
In fact, Achlioptas' proof settles the case $\cR^k$ for any irreducible additive induced-hereditary $\cR$. 

\section{Preliminaries}
We consider only simple finite graphs.
We write $G\le H$ when $G$ is an induced subgraph of $H$.
We identify a graph property with the set of graphs that have that property.
A property $\cP$ is {\em additive}, or {\em (induced-)hereditary}, if it is closed under taking vertex-disjoint unions, or (induced-)subgraphs.
The properties we consider contain the null graph $K_0$ and at least one, but not all (finite simple non-null) graphs.

A {\em $(\cP,\cQ)$-colouring} of $G$ is a partition of $V(G)$ into red and blue vertices, such that the red vertices induce a subgraph $G_\cP \in \cP$, and the blue vertices induce a subgraph $G_\cQ \in \cQ$. 
The {\em product} of $\cP$ and $\cQ$ is $\cP \circ \cQ$, the set of $(\cP,\cQ)$-colourable graphs. We use $(\cP, \cQ)$-colouring, $(\cP, \cQ)$-partition and $(\cP \circ \cQ)$-recognition interchangeably.

Let $\cal P$ be an additive induced-hereditary property.  Then $\cP$ is {\em reducible} if it is the product of two additive induced-hereditary properties; otherwise it is {\em irreducible}. It is true, though by no means obvious, that if $\cP$ is the product of \emph{any} two properties, then it is also the product of two additive induced-hereditary properties~\cite{uft-x}.

The set of minimal forbidden induced-subgraphs for $\cP$ is $\cF(\cP) := \{H \not\in \cP \mid \forall\,G < H,\ G \in \cP\}$. Note that $\cF(\cO) = \{K_2\}$, while all other properties have forbidden subgraphs with at least $3$ vertices. $\cP$ is additive iff every graph in $\cF(\cP)$ is connected. Every hereditary property is induced-hereditary, and the product of additive (induced-hereditary) properties is additive (induced-hereditary).


A graph $H$ is {\em strongly uniquely $(\cP_1, \ldots, \cP_n)$-partitionable} if there is exactly one \emph{ordered} partition $(V_1, \ldots, V_n)$ of $V(H)$ such that for all $i$, $H[V_i] \in \cP_i$. 
More precisely, suppose $V(H) = U_1 \cup \cdots \cup U_n$, where $H[U_i] \in \cP_i$ for all $i$. Then 

\noindent (a) there is a permutation $\phi$ of $\{1, \ldots, n\}$ such that $V_i = U_{\phi(i)}$; 

\noindent (b) if $i, j$ are in the same cycle of $\phi$, then $\cP_i = \cP_j$. 

When the $\cP_i$'s are additive induced-hereditary and irreducible, Mih\'{o}k~\cite{uft-2} gave a construction that can easily be adapted (cf.~\cite[Thm. 5.3]{discussiones}, \cite{uft-x}, \cite{uni-1}) to give a strongly uniquely $(\cP_1, \ldots, \cP_n)$-partitionable graph $H$ with $V_n \not= \emptyset$. 
We use $H$ to show that $\cA \circ \cB$-recognition is at least as hard as 
$\cA$-recognition, when $\cA$ and $\cB$ are additive induced-hereditary properties
(the result is not true for all properties, e.g., $\cB := \{G \mid |V(G)| \geq 10\}$).

\resnum{Theorem\label{complexity-2}}
{Let $\cA$ and $\cB$ be additive induced-hereditary properties.
Then there is a polynomial-time transformation from 
the $\cA$-recognition problem to the $(\cA \circ \cB)$-recognition problem.
}
It is clearly enough to prove this when $\cB$ is irreducible. For any
graph $G$ we will construct (in time linear in $|V(G)|$) a graph $G'$
such that $G \in \cA$ if and only if $G' \in \cA \circ \cB$.

Let $\cA = \cP_1 \circ \cdots \circ \cP_{n-1},\ \cB = \cP_{n}$,
where the $\cP_i$'s are irreducible additive induced-hereditary
properties. Let $H$ be a fixed 
strongly uniquely $(\cP_1, \ldots, \cP_{n})$-partitionable graph, with
partition $(V_1, \ldots, V_{n})$, such that $V_{n} \not= \emptyset$.
Let $v_H$ be some fixed vertex in $V_1$.

For any graph $G$, we construct $G'$ by taking a copy of $G$ and a copy of $H$, and making every vertex of $G$ adjacent to every vertex of $N(v_H) \cap V_{n}$. 
By additivity of $\cA$, if $G$ is in $\cA$, then $G'$ is in $\cA \circ \cB$. 

Conversely, if $G'\in{\cal A}\circ{\cal B}=\cP_1\circ\cdots\circ\cP_{n}$, then it
has an ordered partition $(W_1, \ldots, W_n)$ with $W_i \in \cP_i$ for each $i$.
Since the $\cP_i$'s are induced-hereditary, $G'[W_i]\in\cP_i$ implies $G'[W_i\cap V(H)]\in\cP_i$. Then%
\footnote{Up to some permutation of the subscripts as in (a), (b).}
$(W_1\cap V(H),\dots,W_{n}\cap V(H)) = (V_1, \ldots, V_n)$; in particular, $v_H \in W_1$.

Suppose some $w \in V(G)$ is in $W_n$. 
Now $(V_1 \setminus \{v_H\}, V_2, \ldots, V_{n-1}, V_n \cup \{w\})$ is a $(\cP_1,\dots,\cP_{n})$-partition of $(H-v_H)+w \cong H$. Then $(V_1 \setminus \{v_H\}, V_2, \ldots, V_{n-1},$ $V_n \cup \{v_H\})$ is a 
$(\cP_1,\dots,\cP_{n})$-partition of $H$ that is different from $(V_1, \ldots, V_n)$ (since $V_n \not= \emptyset$), a contradiction.

Thus no vertex of $G$ is in $W_n$, and so $G \leq G'[W_1 \cup \cdots \cup W_{n-1}] \in 
\cP_1\circ\cdots\circ\cP_{n-1}=\cA$, and $G \in \cA$ as required.
\resnumend
\newline

We will prove the main result by transforming {$p$-\sc{in}-$r$-\sc{SAT}} to $(\cP,\cQ)$-colouring, where $p$ and $r$ are fixed integers depending on $\cP$ and $\cQ$. 
Schaefer~\cite{schaefer} showed {$p$-\sc{in}-$r$-\sc{SAT}} to be NP-complete, even for formulae with all literals unnegated, for any fixed $p$ and $r$, so long as $1 \leq p < r$ and $r \geq 3$. We restate it as:
\newline

\noindent {$p$-\sc{in}-$r$-\sc{colouring}}

\noindent {\bf Instance}: an $r$-uniform hypergraph.

\noindent {\bf Problem}: is there a set of vertices $U$ such that, for each hyper-edge $e$, $|U \cap e| = p$?

%
%
%

\section{NP-hardness\label{sec-nphard}}
\begin{figure}[htb]
\begin{center}
\input{A-gadget.pstex_t} 
\caption{The forbidden graphs $F_{\cP}$ and $F_{\cQ}$, and the replicator gadget $R$.}
\label{Fig-same-gadget}
\end{center}
\end{figure}

\resnum{Theorem\label{complexity-3}}
{Let $\cP$ and $\cQ$ be additive induced-hereditary properties, $\cP \circ \cQ \not= \cO^2$. Then $(\cP \circ \cQ)$-recognition is NP-hard. Moreover, it is NP-complete iff $\cP$- and $\cQ$-recognition are both in NP.
}
We will prove the first part; the second part then follows by 
Theorem~\ref{complexity-2}. Also by Theorem~\ref{complexity-2}
(and by the well-known NP-hardness of recognising 
$\cO^3$~\cite{karp}\label{karp-1}), 
we need only consider the case where $\cP$ and $\cQ$ are irreducible.
By Theorem~\ref{complexity-2} there is a strongly uniquely $(\cP, \cQ)$-colourable graph $G_{\cP, \cQ}$ that we use to ``force'' vertices to be in $\cP$ or $\cQ$. 

More formally, let the unique partition be $V(G_{\cP, \cQ}) = U_\cP \cup U_\cQ$. Choose $p \in U_\cP$. If $G_{\cP, \cQ} \leq H$, and $v \not\in V(G_{\cP, \cQ})$ satisfies $N(v) \cap U_\cQ = N(p) \cap U_\cQ$, then in any 
$(\cP, \cQ)$-colouring of $H$, $v$ must be in the $\cP$-part%
\footnote{To be precise, we mean that $v$ is coloured the same as $p$: if $\cP = \cQ$ then a $(\cP, \cQ)$-colouring is also a $(\cQ, \cP)$-colouring, but we adopt the convention that the $\cP$-part is the part containing $p$.}; 
otherwise, in $G_{\cP, \cQ}$ we could transfer $p$ over to the $\cQ$ part, giving us a different $(\cP, \cQ)$-colouring. Similarly we choose $q \in U_\cQ$, whose neighbours we use to force vertices to be in $\cQ$. $G_{\cP, \cQ}$ is our first gadget.
\newline

An {\em end-block} of a graph $G$ is a block of $G$ that contains at most one 
cut-vertex of $G$; if $G$ has no cut-vertices, then $G$ is itself an end-block. 
Let $B_\cP$ be an end-block of $F_\cP \in \cF(\cP)$, chosen to have the least number of vertices among all the end-blocks of all the graphs in $\cF(\cP)$
(see Figure~\ref{Fig-same-gadget}). 
Because $\cP$ is additive and non-trivial, $F_\cP$ is connected and has at least two vertices, so $B_\cP$ has $k \geq 2$ vertices. The point to note is that, if $H$ is a graph in $\cP$, then adding an end-block with fewer than $k$ vertices produces another graph in $\cP$.

Let $y_\cP$ be the unique cut-vertex contained in $B_\cP$ (if $B_\cP = F_\cP$, pick $y_\cP$ arbitrarily), and let $x_\cP$ be a vertex of $B_\cP$ adjacent to $y_\cP$. 
Let $F'_\cP$ be the graph obtained by adding an extra copy of $B_\cP$ (incident to the same cut-vertex $y_\cP$), and let $x'_\cP$ be a vertex in this new copy that is adjacent to $y_\cP$.

Similarly, we choose $B_\cQ$ to be an end-block of $F_\cQ \in \cF(\cQ)$, minimal among the end-blocks of 
graphs in $\cF(\cQ)$; we add a copy of $B_\cQ$, and pick $x_\cQ$, $y_\cQ$ and $x'_\cQ$ as above. We identify $x_\cP$ with $x_\cQ$, $y_\cP$ with $y_\cQ$, $x'_\cP$ with $x'_\cQ$, and label the identified vertices $x, y, x'$. 

Finally, we force all the vertices of $F'_\cP$ (except for $x, y, x'$) to be in $\cP$, and all the vertices of $F'_\cQ$ (except for $x, y, x'$) to be in $\cQ$. That is, we add a copy of $G_{\cP, \cQ}$, and make every vertex of $F'_{\cal P}-\{x,y,x'\}$ adjacent to every vertex of $N(p) \cap U_\cQ$, and every vertex of $F'_{\cal Q}-\{x,y,x'\}$ adjacent to every vertex of $N(q) \cap U_\cP$ (cf. Figure~\ref{Fig-same-gadget}).


It can be 
readily checked that the resulting gadget $R$ (for `replicator') has the following properties:
\begin{itemize}
\item In a $(\cP, \cQ)$-colouring of $R$, if $x$ is in $\cP$, then $y$ is in $\cQ$ and $x'$ is in $\cP$; similarly, if $x$ is in $\cQ$, then $y$ is in $\cP$ and $x'$ is in $\cQ$. 
So $x$ and $x'$ always have the same colour, that is different from that of $y$. Moreover, there is at least one  colouring (in fact, exactly one) in which $x$ and $x'$ are in $\cP$, and at least one in which both are in $\cQ$.
\item Identify $x$ with a vertex $z$ of some graph $H$ to obtain $H_R$; then 
$(\cP, \cQ)$-colourings of $H$ and $R$ that agree on $x$, together give a 
$(\cP, \cQ)$-colouring of $H_R$. We can then similarly identify $x'$ with some vertex $z'$ of a graph $H'$, and attach more copies of $R$ at $x$ or $x'$.
\end{itemize}
\clearpage

We thus have a gadget that ``replicates'' the colour of $x$ on $x'$, while preserving valid colourings.
\newline

Let $H_\cP$ be a forbidden subgraph for $\cP$ with the least possible number of vertices, say $p+1$; similarly choose $H_\cQ \in \cF(\cQ)$ on $q+1$ vertices, where $q+1$ is as small as possible, so any graph on at most $p$ (resp. $q$) vertices is in $\cP$ (resp. $\cQ$). Since $\cP$ and $\cQ$ are not both $\cO$, $p+q \geq 3$, and so {\sc $p$-in-$(p+q)$-colouring} is NP-complete. We will construct a third gadget to transform this to $(\cP, \cQ)$-colouring.

We start with an independent set $S$ on $p+q$ vertices, $\{x_1, \ldots, x_{p+q}\}$. For every $(p+1)$-subset of $S$, say $T_j = \{x_1, \ldots, x_{p+1}\}$, add a 
disjoint 
copy of $H_\cP$ whose vertices are labeled $x^j_1, \ldots, x^j_{p+1}$. For each $i = 1, \ldots, p+1,$ use a new copy $R_{i,j}$ of $R$ to ensure that $x_i$ and $x^j_i$ are always coloured the same; to do this, identify the vertices $x$ and $x'$ of $R_{i,j}$ with $x_i$ and $x^j_i$. For every $(q+1)$-subset of $S$ we add a copy of $H_\cQ$ in the same manner. 
Thus every vertex $x_i \in S$ will have $\ell = {p+q-1 \choose p} + {p+q-1 \choose q}$ 
`shadow vertices' $x_i^1, \ldots, x_i^{\ell}$ from copies of $H_\cP$ and $H_\cQ$. Call this gadget $N$ (for `pin cushion' --- the copies of $H_\cP$ and $H_\cQ$ being stuck into the independent set $S$ by `pins' or `replicators').

In a $(\cP, \cQ)$-colouring of $N$, no $p+1$ vertices of $S$ can be in $\cP$, and no $q+1$ vertices can be in $\cQ$, so exactly $p$ vertices of $S$ are in $\cP$, and exactly $q$ are in $\cQ$. Conversely, suppose that exactly $p$ vertices of $S$ are coloured red, and  the other $q$ are blue; colour each vertex $x^j_i$ the same as $x_i$, $1 \leq i \leq  p+q$, $1 \leq j \leq \ell$. Then each copy of $H_\cP$ has at most $p$ red and at most $q$ blue vertices, giving it a valid $(\cP, \cQ)$-colouring. The colouring on the rest of each gadget $R_{i,j}$ is then forced, and we have a $(\cP, \cQ)$-colouring of all of $N$.
\newline

Now, given a $(p+q)$-uniform hypergraph $\cal H$, we stick a copy of $N$ onto every hyper-edge. The resulting graph is $(\cP, \cQ)$-colourable iff $\cal H$ has a {\sc $p$-in-$(p+q)$-colouring}.
\resnumend
\newline

\section{New directions}
How far can the main result be extended? Uniquely 
$(\cP_1, \ldots, \cP_n)$-partitionable graphs exist even in many cases where the $\cP_i$'s are not additive~\cite{corr}; however, this includes finite $\cP_i$'s, so the existence of uniquely colourable graphs does not guarantee NP-hardness.

It may be useful to restate the result as follows: if the graphs in $\cF(\cP)$ and $\cF(\cQ)$ are all connected, then $(\cP, \cQ)$-colouring is NP-hard. This is also true if the graphs in $\cF(\cP)$ and $\cF(\cQ)$ are all disconnected, since $G \in \cP \circ \cQ \Leftrightarrow \overline{G} \in \overline{\cP} \circ \overline{\cQ}$, where $\overline{\cP}$ is defined by $\cF(\overline{\cP}) := \{\overline{H} \mid H \in \cF(\cP)\}$.

A natural problem to tackle next would be classifying the complexity of $\cR^k$-recognition, where $\cR$ has both connected and disconnected minimal forbidden induced-subgraphs. One of the simplest such cases is $\cR = (\cO \cup \cK)$, where $\cK$ is the set of all cliques: $\cF(\cO \cup \cK) = \{P_3, \overline{P_3}\}$. Gimbel et al.~\cite{gimbel} noted that $G \in \cO^k \Leftrightarrow nG \in (\cO \cup \cK)^k$ (where $n = |V(G)|$); so $(\cO \cup \cK)^k$-recognition is NP-complete for $k \geq 3$ (and, in fact, polynomial for $k = 1,2$). 

Another natural problem is $(\cP, \cQ)$-colouring, where all graphs in $\cF(\cP)$ are connected, and all those in $\cF(\cQ)$ are disconnected. 
In all problems, it may make sense to restrict attention to hereditary properties with finitely many forbidden subgraphs.

Another class of problems often considered in the literature is $(\cD: \cP)$-recognition: given a graph $G$ in the domain $\cD$, is $G$ in $\cP$? This is just $(\cD \cap \cP)$-recognition; if $\cD$ and $\cP$ are both additive induced-hereditary, then so is $\cD \cap \cP$, with $\cF(\cD \cap \cP) = \min_{\leq}(\cF(\cD) \cup \cF(\cP))$. We leave it as an open question, for reducible $\cP$, to determine when $\cD \cap \cP$ is also reducible; Mih\'ok's characterisations~\cite{uft-1, uft-2} of reducibility may prove useful.

\section{Notes and acknowledgements}
The most important part of the proof is the `replicator' gadget. Phelps and R\"odl~\cite[Thm. 6.2]{phelps} and Brown~\cite[Thm. 2.3]{brown} used different gadgets to perform similar roles. The forcing technique of Theorem~\ref{complexity-2} was first used in~\cite[Thm. 2]{comp3} and~\cite[Lemma 3]{uni-1}.

Contacts with Lozin were very helpful, as they spurred the author to look at ($K_m$-free, $K_n$-free)-colouring, not knowing it had been settled in~\cite{cai}. Kratochv\'{\i}l and Schiermeyer~\cite{comp3} proved a special case of Theorem~\ref{complexity-3} that covered the case $m=2$;
($K_2$-free, $K_n$-free)-colouring; 
I started my proof for general $m$ and $n$ by adapting theirs, and ended up strengthening and simplifying it considerably.

I would like to thank Bruce Richter for many helpful conversations, detailed comments that improved the presentation of the paper, and for spotting a flaw in my original `pin cushion' gadget. 
The result here forms part of the Ph.D.\ thesis that I am writing under his supervision. 
I would also like to thank the Canadian government for fully funding my studies through a Commonwealth Scholarship.





\end{document}

%% file: A-gadget.pstex_t
\begin{picture}(0,0)%
\epsfig{file=A-gadget.pstex}%
\end{picture}%
\setlength{\unitlength}{3947sp}%
\begingroup\makeatletter\ifx\SetFigFont\undefined%
\gdef\SetFigFont#1#2#3#4#5{%
  \reset@font\fontsize{#1}{#2pt}%
  \fontfamily{#3}\fontseries{#4}\fontshape{#5}%
  \selectfont}%
\fi\endgroup%
\begin{picture}(5408,8409)(143,-7573)
\put(2701,-4211){\makebox(0,0)[lb]{\smash{\SetFigFont{12}{14.4}{\rmdefault}{\mddefault}{\updefault}$x'$}}}
\put(2451,-7561){\makebox(0,0)[lb]{\smash{\SetFigFont{12}{14.4}{\rmdefault}{\mddefault}{\updefault}$G_{\cP, \cQ}$}}}
\put(2863,-4811){\makebox(0,0)[lb]{\smash{\SetFigFont{12}{14.4}{\rmdefault}{\mddefault}{\updefault}$y$}}}
\put(2201,-7261){\makebox(0,0)[lb]{\smash{\SetFigFont{12}{14.4}{\rmdefault}{\mddefault}{\updefault}$q$}}}
\put(5551,-5911){\makebox(0,0)[lb]{\smash{\SetFigFont{12}{14.4}{\rmdefault}{\mddefault}{\updefault}$R$}}}
\put(3201,-944){\makebox(0,0)[lb]{\smash{\SetFigFont{12}{14.4}{\rmdefault}{\mddefault}{\updefault}$B_{\cQ}$}}}
\put(1051,-1501){\makebox(0,0)[lb]{\smash{\SetFigFont{12}{14.4}{\rmdefault}{\mddefault}{\updefault}$F_{\cP}$}}}
\put(1601,-961){\makebox(0,0)[lb]{\smash{\SetFigFont{12}{14.4}{\rmdefault}{\mddefault}{\updefault}$B_{\cP}$}}}
\put(4301,-1501){\makebox(0,0)[lb]{\smash{\SetFigFont{12}{14.4}{\rmdefault}{\mddefault}{\updefault}$F_{\cQ}$}}}
\put(3051,-61){\makebox(0,0)[lb]{\smash{\SetFigFont{12}{14.4}{\rmdefault}{\mddefault}{\updefault}$y_\cQ$}}}
\put(3051,-514){\makebox(0,0)[lb]{\smash{\SetFigFont{12}{14.4}{\rmdefault}{\mddefault}{\updefault}$x_\cQ$}}}
\put(1951,-61){\makebox(0,0)[lb]{\smash{\SetFigFont{12}{14.4}{\rmdefault}{\mddefault}{\updefault}$y_\cP$}}}
\put(1951,-511){\makebox(0,0)[lb]{\smash{\SetFigFont{12}{14.4}{\rmdefault}{\mddefault}{\updefault}$x_\cP$}}}
\put(3043,-6141){\makebox(0,0)[lb]{\smash{\SetFigFont{12}{14.4}{\rmdefault}{\mddefault}{\updefault}$p$}}}
\put(2701,-5351){\makebox(0,0)[lb]{\smash{\SetFigFont{12}{14.4}{\rmdefault}{\mddefault}{\updefault}$x$}}}
\end{picture}

%% file: nphard.bbl
\begin{thebibliography}{24}
\bibitem{comp1} D.\ Achlioptas, The complexity of G-free
colourability, {\em Discrete Math.\ } {\bf 165-166} (1997) 21--30.\ 

\bibitem{subchrom} M.O.\ Albertson, R.E.\ Jamison, S.T.\ Hedetniemi, S.C.\ Locke, The subchromatic number of a graph, {\em Discrete Math} {\bf 74} (1989) 33--49.


\bibitem{berger} A.J.\  Berger, Minimal forbidden subgraphs of reducible graph properties, {\em Discuss.\ Math.\ Graph Theory} {\bf 21} (2001) 111-117.

\bibitem{brand} A.\ Brandst\"adt, V.B.\ Le, T.\ Szymcak, The complexity of some problems related to {\sc graph $3$-colorability}, {\em Disc.\ Appl.\ Math.\ } {\bf 89} (1998) 59--73.

\bibitem{uni-1} I.\ Broere and J.\ Bucko, Divisibility in additive
hereditary properties and uniquely partitionable graphs, {\em Tatra
Mt.\ Math.\ Publ.\ } {\bf 18} (1999) 79--87.\ 

\bibitem{brown-thesis} J.I.\ Brown, A theory of generalized graph colourings, Ph.\ D.\ Thesis, Department of Mathematics, University of Toronto (1987).

\bibitem{brown} J.I.\ Brown, The complexity of generalized graph
colorings,  {\em Discrete Appl.\ Math.\ } {\bf 69} (1996)  257--270.

\bibitem{brown-perfect} J.I.\ Brown and D.G.\ Corneil, Perfect colourings, {\em Ars Combin.\ }{\bf 30} (1990) 141--159.

\bibitem{cai} L.\ Cai and D.G.\ Corneil, A generalization of perfect 
graphs---$i$-perfect graphs, {\em J.\ Graph Theory} {\bf 23} (1996) 87--103.

\bibitem{discussiones} A.\ Farrugia and R.B.\ Richter, Unique factorisation of additive induced-hereditary properties, to appear in {\em Discuss.\ Math.\ Graph Theory}.

\bibitem{uft-x} A.\ Farrugia and R.B.\ Richter, Factorisation, reducibility, 
co-primality, and uniquely colourable graphs, in preparation.

\bibitem{corr} A.\ Farrugia and R.B.\ Richter, Unique factorisation
of induced-hereditary disjoint compositive properties, Research Report CORR 
2002-ZZ (2002) Department of Combinatorics and Optimization,

\bibitem{subchrom-2} J.\ Fiala, K.\ Jansen, V.B.\ Le and E.\ Seidel, Graph subcolorings: complexity and algorithms, {\em Lecture Notes in Computer Science} {\bf 2204} (Proceedings, Boltenhagen, 2001) 154--165.


\bibitem{g-j-stock} M.R.\ Garey, D.S.\ Johnson and L.\ Stockmeyer, Some simplified NP-complete problems, {\em Theor.\ Comput.\ Sci.\ } {\bf 1} (1976) 237--267.

\bibitem{gimbel} J.\ Gimbel, D.\ Kratsch and L.\ Stewart, On cocolourings and cochromatic numbers of graphs, {\em Disc.\ Appl.\ Math.} {\bf 48} (1994) 111--127.

\bibitem{arbor} S.L.\ Hakimi and E.F.\ Schmeichel, A note on the vertex arboricity of a graph, {\em SIAM J.\ Discrete Math.} {\bf 2} (1989) 64--67.

\bibitem{karp} R.M.\ Karp, Reducibility among combinatorial problems, in R.E.\ Miller and J.W.\ Thatcher (eds.), {\em Complexity of computer computations}, Plenum Press, New York, 85--103.

\bibitem{comp3} J.\ Kratochv\'{\i}l and I.\ Schiermeyer, On the
computational complexity of $(O,P)$-partition problems, {\em Discuss.\
Math.\ Graph Theory} {\bf 17} (1997) 253--258.

\bibitem{uft-1} P.\ Mih\'{o}k, G.\ Semani\v{s}in and R.\ Vasky, Additive and hereditary properties of graphs are uniquely factorizable into irreducible factors, {\em J.\ Graph Theory} {\bf 33} (2000) 44--53.\ 

\bibitem{uft-2} P.\ Mih\'{o}k, Unique Factorization Theorem,
{\em Discuss.\ Math.\ Graph Theory} {\bf 20} (2000) 143--153.\  

\bibitem{monien} B.\ Monien, correspondence with Brandst\"adt, Le and Szymcak, 1984.

\bibitem{phelps} K.T.\ Phelps and V.\ R\"odl, 
Algorithmic complexity of coloring simple hypergraphs and Steiner triple systems, {\em Combinatorica} {\bf 4} (1984) 79--88.

\bibitem{schaefer} T.J.\ Schaefer, The complexity of satisfiability problems, {\em Proc.\ 10th Ann.\ ACM Symp.\ on Theory of Computing}, Association for Computing Machinery, New York (1978) 216--226.

\bibitem{szi-t} J.\ Szigeti and Z.\ Tuza, Generalized colorings and avoidable orientations, {\em Discuss.\ Math.\ Graph Theory} {\bf 17} (1997) 137--146.

\end{thebibliography}
